\documentclass[12pt,fleqn,leqno]{article}

\author{
D.D.\ Poro\c sniuc\\  }
\date{}
\title{A locally symmetric K\"ahler Einstein structure on the cotangent
bundle of a space form}

\begin{document}

\maketitle
\begin{abstract}
We obtain a locally symmetric K\"ahler Einstein structure on the
cotangent bundle of a Riemannian manifold of negative constant
sectional curvature. Similar results are obtained on a tube around
zero section in the cotangent bundle, in the case of a Riemannian
manifold of  positive constant sectional curvature. The obtained
K\"ahler Einstein structures cannot have constant holomorphic
sectional curvature.


\end{abstract}

\vskip5mm {\large \bf Introduction} \vskip5mm

 The differential geometry of the cotangent
bundle $T^*M$ of a Riemannian manifold $(M,g)$ is almost similar
to that of the tangent bundle $TM$. However there are some
differences, because the lifts (vertical, complete, horizontal
etc.) to $T^*M$ cannot be defined just like in the case of $TM$.

In \cite{OprPor} V. Oproiu and the present author have obtained a
family of natural K\"ahler Einstein structures $(G,J)$ of diagonal
type induced on $T^*M$ from the Riemannian metric $g$. The
obtained K\"ahler structure on $T^*M$ depends on one essential
parameter $u$, which is a smooth function depending on the energy
density $t$ on $T^*M$. If the K\"ahler structure is Einstein they
get a second order differential equation fulfilled by the
parameter $u$. In the case of the general solution, they have
obtained that $(T^*M,G,J)$ has constant holomorphic sectional
curvature.

 In this paper we study the singular case where the parameter
$u=A,A\in {\bf R}.$ The considered natural Riemannian metric $G$
of diagonal type on $T^*M$ is defined by using one parameter $v$
which is a smooth function depending on the energy density $t$ on
$T^*M$. The vertical distribution $VT^*M$ and the horizontal
distribution $HT^*M$ are orthogonal to each other but the dot
products induced on them from $G$ are not isomorphic (isometric).

 Next, the natural almost complex structures $J$ on
$T^*M$ that interchange the vertical and horizontal distributions
depends of one essential parameter $v$.

 After that, we obtain that $G$ is Hermitian with respect to $J$ and it
 follows that the fundamental $2$-form $\phi$, associated to the
 almost Hermitian structure $(G,J)$ is the fundamental form defining the
 usual symplectic structure on $T^*M$, hence it is closed.

  From the integrability condition for $J$ it
follows that the base manifold $M$ must have constant sectional
curvature $c$ and the parameter $v$ must be constant.

 If the constant
sectional curvature c is negative then we obtain a locally
symmetric K\"ahler Einstein structure defined on the whole $T^*M$.
If the constant sectional curvature c of $M$ is positive then we
get a similar structure defined on a tube around zero section in
$T^*M$. The K\"ahler Einstein manifolds obtained cannot have
constant holomorphic sectional curvature.

 The manifolds, tensor fields and geometric objects we consider
in this paper, are assumed to be differentiable of class
$C^{\infty}$ (i.e. smooth). We use the computations in local
coordinates but many results from this paper may be expressed in
an invariant form. The well known summation convention is used
throughout this paper, the range for the indices $h,i,j,k,l,r,s$
being always${\{}1,...,n{\}}$ (see \cite{GheOpr}, \cite{OprPap1},
\cite{OprPap2}, \cite{YanoIsh}). We shall denote by
${\Gamma}(T^*M)$ the module of smooth vector fields on $T^*M$.

\vskip5mm {\large \bf 1. Some geometric properties of $T^*M$}
\vskip5mm

Let $(M,g)$ be a smooth $n$-dimensional Riemannian manifold and
denote its cotangent bundle by $\pi :T^*M\longrightarrow M$.
Recall that there is a structure of a $2n$-dimensional smooth
manifold on $T^*M$, induced from the structure of smooth
$n$-dimensional manifold  of $M$. From every local chart
$(U,\varphi )=(U,x^1,\dots ,x^n)$  on $M$, it is induced a local
chart $(\pi^{-1}(U),\Phi )=(\pi^{-1}(U),q^1,\dots , q^n,$
$p_1,\dots ,p_n)$, on $T^*M$, as follows. For a cotangent vector
$p\in \pi^{-1}(U)\subset T^*M$, the first $n$ local coordinates
$q^1,\dots ,q^n$ are  the local coordinates $x^1,\dots ,x^n$ of
its base point $x=\pi (p)$ in the local chart $(U,\varphi )$ (in
fact we have $q^i=\pi ^* x^i=x^i\circ \pi, \ i=1,\dots n)$. The
last $n$ local coordinates $p_1,\dots ,p_n$ of $p\in \pi^{-1}(U)$
are the vector space coordinates of $p$ with respect to the
natural basis $(dx^1_{\pi(p)},\dots , dx^n_{\pi(p)})$, defined by
the local chart $(U,\varphi )$,\ i.e. $p=p_idx^i_{\pi(p)}$.

An $M$-tensor field of type $(r,s)$ on $T^*M$ is defined by sets
of $n^{r+s}$ components (functions depending on $q^i$ and $p_i$),
with $r$ upper indices and $s$ lower indices, assigned to induced
local charts $(\pi^{-1}(U),\Phi )$ on $T^*M$, such that the local
coordinate change rule is that of the local coordinate components
of a tensor field of type $(r,s)$ on the base manifold $M$ (see
\cite{Mok} for further details in the case of the tangent bundle).
An usual tensor field of type $(r,s)$ on $M$ may be thought of as
an $M$-tensor field of type $(r,s)$ on $T^*M$. If the considered
tensor field on $M$ is covariant only, the corresponding
$M$-tensor field on $T^*M$ may be identified with the induced
(pullback by $\pi $) tensor field on $T^*M$.

Some useful $M$-tensor fields on $T^*M$ may be obtained as
follows. Let $v,w:[0,\infty ) \longrightarrow {\bf R}$ be a smooth
functions and let $\|p\|^2=g^{-1}_{\pi(p)}(p,p)$ be the square of
the norm of the cotangent vector $p\in \pi^{-1}(U)$ ($g^{-1}$ is
the tensor field of type (2,0) having the components $(g^{kl}(x))$
which are the entries of the inverse of the matrix $(g_{ij}(x))$
defined by the components of $g$ in the local chart $(U,\varphi
)$). The components $g_{ij}(\pi(p))$, $p_i$, $v(\|p\|^2)p_ip_j $
define $M$-tensor fields of types $(0,2)$, $(0,1)$, $(0,2)$ on
$T^*M$, respectively. Similarly, the components $g^{kl}(\pi(p))$,
$g^{0i}=p_hg^{hi}$, $w(\|p\|^2)g^{0k}g^{0l}$ define $M$-tensor
fields of type $(2,0)$, $(1,0)$, $(2,0)$ on $T^*M$, respectively.
Of course, all the components considered above are in the induced
local chart $(\pi^{-1}(U),\Phi)$.

 The Levi Civita connection $\dot \nabla $ of $g$ defines a direct
sum decomposition

\begin{equation}
TT^*M=VT^*M\oplus HT^*M.
\end{equation}
of the tangent bundle to $T^*M$ into vertical distributions
$VT^*M= {\rm Ker}\ \pi _*$ and the horizontal distribution
$HT^*M$.

 If $(\pi^{-1}(U),\Phi)=(\pi^{-1}(U),q^1,\dots ,q^n,p_1,\dots ,p_n)$
is a local chart on $T^*M$, induced from the local chart
$(U,\varphi )= (U,x^1,\dots ,x^n)$, the local vector fields
$\frac{\partial}{\partial p_1}, \dots , \frac{\partial}{\partial
p_n}$ on $\pi^{-1}(U)$ define a local frame for $VT^*M$ over $\pi
^{-1}(U)$ and the local vector fields $\frac{\delta}{\delta
q^1},\dots ,\frac{\delta}{\delta q^n}$ define a local frame for
$HT^*M$ over $\pi^{-1}(U)$, where
$$
\frac{\delta}{\delta q^i}=\frac{\partial}{\partial
q^i}+\Gamma^0_{ih} \frac{\partial}{\partial p_h},\ \ \ \Gamma
^0_{ih}=p_k\Gamma ^k_{ih}
 $$
and $\Gamma ^k_{ih}(\pi(p))$ are the Christoffel symbols of $g$.

The set of vector fields $(\frac{\partial}{\partial p_1},\dots
,\frac{\partial}{\partial p_n}, \frac{\delta}{\delta q^1},\dots
,\frac{\delta}{\delta q^n})$ defines a local frame on $T^*M$,
adapted to the direct sum decomposition (1).

We consider
\begin{equation}
t=\frac{1}{2}\|p\|^2=\frac{1}{2}g^{-1}_{\pi(p)}(p,p)=\frac{1}{2}g^{ik}(x)p_ip_k,
\ \ \ p\in \pi^{-1}(U)
\end{equation}
the energy density defined by $g$ in the cotangent vector $p$. We
have $t\in [0,\infty)$ for all $p\in T^*M$.

From now on we shall work in a fixed local chart $(U,\varphi)$ on
$M$ and in the induced local chart $(\pi^{-1}(U),\Phi)$ on $T^*M$.

\vskip5mm {\large \bf 2. An almost K\"ahler structure on the
$T^*M$} \vskip5mm

Consider a real valued smooth function $v$ defined on
$[0,\infty)\subset {\bf R}$ and a real constant A. We define the
following symmetric $M$-tensor field of type $(0,2)$  on $T^*M$
having the components
\begin{equation}
G_{ij}(p)= A g_{ij}(\pi(p))+v(t)p_ip_j.
\end{equation}

It follows easily that the matrix $(G_{ij})$ is positive definite
if and only if $A>0,\ A+2tv>0$. The inverse of this matrix has the
entries

\begin{equation}
H^{kl}(p)= \frac{1}{A} g^{kl}(\pi(p))+w(t)g^{0k}g^{0l},
\end{equation}
where
\begin{equation}
w= -\ \frac{v}{A(A+2tv)}.
\end{equation}

The components $H^{kl}$ define an symmetric $M$-tensor field of
type $(2,0)$
on $T^*M$.\\

{\bf Remark.} If the matrix $(G_{ij})$ is positive definite, then
its inverse $(H^{kl})$ is positive definite too.

Using the $M$-tensor fields defined by $G_{ij},\ H^{kl}$, the
following Riemannian metric may be considered on $T^*M$
\begin{equation}
G=G_{ij}dq^idq^j+H^{ij}Dp_iDp_j,
\end{equation}
where $Dp_i=dp_i-\Gamma^0_{ij}dq^j$ is the absolute (covariant)
differential of $p_i$ with respect to the Levi Civita connection
$\dot\nabla$ of $g$ . Equivalently, we have
$$
G(\frac{\delta}{\delta q^i},\frac{\delta}{\delta
q^j})=G_{ij},~~~G(\frac{\partial}{\partial p_i}
,\frac{\partial}{\partial p_j})=H^{ij},~~
G(\frac{\partial}{\partial p_i},\frac{\delta}{\delta q^j})=~
G(\frac{\delta}{\delta q^j},\frac{\partial}{\partial p_i})=0.
$$

Remark that $HT^*M,~VT^*M$ are orthogonal to each other with
respect to $G$, but the Riemannian metrics induced from $G$ on
$HT^*M,~VT^*M$ are not the same, so the considered metric $G$ on
$T^*M$ is not a metric of Sasaki type.
 Remark also that the system of 1-forms
$(dq^1,...,dq^n,Dp_1,...,Dp_n)$ defines a local frame on
$T^{*}T^*M$, dual to the local frame $(\frac{\delta}{\delta
q^1},\dots ,\frac{\delta}{\delta q^n},\frac{\partial}{\partial
p_1},\dots ,\frac{\partial}{\partial p_n})$ adapted to the direct
sum decomposition (1).

Next, an almost complex structure $J$ is defined on $T^*M$ by the
same $M$-tensor fields $G_{ij},\ H^{kl}$, expressed in the adapted
local frame by
\begin{equation}
J\frac{\delta}{\delta
q^i}=G_{ik}\frac{\partial}{\partial p_k},\ \ \
J\frac{\partial}{\partial p_i}=-H^{ik}\frac{\delta}{\delta q^k}
\end{equation}

 From the property
of the $M$-tensor field $H^{kl}$ to be defined by the inverse of
the matrix defined by the components of the $M$-tensor field
$G_{ij}$, it follows easily that $J$ is an almost complex
structure on $T^*M$. \\

{\bf Theorem 1}. \it $(T^*M,G,J)$ is an almost K\"ahler
manifold.\\

 Proof. \rm Since the matrix $(H^{kl})$ is the inverse of the
 matrix $(G_{ij})$, it follows easily that
$$
G(J\frac{\delta}{\delta q^i},J\frac{\delta}{\delta
q^j})=G(\frac{\delta}{\delta q^i},\frac{\delta}{\delta q^j}),\ \ \
G(J\frac{\partial}{\partial p_i},J\frac{\partial}{\partial
p_j})=G(\frac{\partial}{\partial p_i},\frac{\partial}{\partial
p_j}),$$ $$ G(J\frac{\partial}{\partial p_i},J\frac{\delta}{\delta
q^j})=G(\frac{\partial}{\partial p_i},\frac{\delta}{\delta
q^j})=0.
$$

Hence
$$G(JX,JY)=G(X,Y),\ \ \forall\ X,Y{\in}{\Gamma}(T^*M).$$
Thus  $(T^*M,G,J)$ is an almost Hermitian manifold.

The fundamental $2$-form associated with this almost Hermitian
structure is $\phi$, defined by
$$\phi(X,Y) = G(X,JY),\ \ \ \forall\ X,Y{\in}{\Gamma}(T^*M).$$
By a straightforward computation we get
$$
\phi(\frac{\delta}{\delta q^i},\frac{\delta}{\delta q^j})=0,\ \ \
\phi(\frac{\partial}{\partial p_i},\frac{\partial}{\partial
p_j})=0,\ \ \ \phi(\frac{\partial}{\partial
p_i},\frac{\delta}{\delta q^j})= \delta^i_j.
$$
Hence
\begin{equation}
\phi =Dp_i\wedge dq^i= dp_i\wedge dq^i,
\end{equation}
due to the symmetry of $\Gamma^0_{ij}=p_h\Gamma^h_{ij}$. It
follows that $\phi$ does coincide with the fundamental $2$-form
defining the usual symplectic structure on $T^*M$. Of course, we
have $d\phi =0$, i.e. $\phi$ is closed. Therefore $(T^*M,G,J)$ is
an almost K\"ahler manifold.

\vskip5mm {\large \bf 3. A K\"ahler structure on $T^*M$} \vskip5mm

We shall study the integrability of the almost complex structure
defined by $J$ on $T^*M$. To do this we need the following well
known formulas for the brackets of the vector fields
$\frac{\partial}{\partial p_i},\frac{\delta}{\delta q^i},~
i=1,...,n$
\begin{equation}
[\frac{\partial}{\partial p_i},\frac{\partial}{\partial
p_j}]=0,~~~[\frac{\partial}{\partial p_i},\frac{\delta}{\delta
q^j}]=\Gamma^i_{jk}\frac{\partial}{\partial p_k},~~~
[\frac{\delta}{\delta q^i},\frac{\delta}{\delta q^j}]
=R^0_{kij}\frac{\partial}{\partial p_k},
\end{equation}
where $R^h_{kij}(\pi(p))$ are the local coordinate components of
the curvature tensor field of $\dot \nabla$ on $M$ and
$R^0_{kij}(p)=p_hR^h_{kij}$ . Of course, the components
 $R^0_{kij}$, $R^h_{kij}$ define M-tensor fields of types
 (0,3), (1,3) on $T^*M$, respectively.\\

{\bf Theorem 2. } {\it The Nijenhuis tensor field of the almost
complex structure $J$ on $T^*M$ is given by}
\begin{equation}
\left\{
\begin{array}{l}
N(\frac{\delta}{\delta q^i},\frac{\delta}{\delta
q^j})=-\{Av(\delta^h_ig_{jk}-
\delta^h_jg_{ik})+R^h_{kij}\}p_h\frac{\partial}{\partial p_k},
\\ \mbox{ } \\
N(\frac{\delta}{\delta q^i},\frac{\partial}{\partial
p_j})=-H^{kl}H^{jr}\{Av(\delta^h_ig_{rl}-
\delta^h_rg_{il})+R^h_{lir}\}p_h\frac{\delta}{\delta q^k},
\\ \mbox{ } \\
N(\frac{\partial}{\partial p_i},\frac{\partial}{\partial
p_j})=-H^{ir}H^{jl}\{Av(\delta^h_lg_{rk}-
\delta^h_rg_{lk})+R^h_{klr}\}p_h\frac{\partial}{\partial p_k}.
\end{array}
\right.
\end{equation}

\vskip2mm

{\it Proof. } Recall that the Nijenhuis tensor field $N$ defined
by $J$ is given by
$$
N(X,Y)=[JX,JY]-J[JX,Y]-J[X,JY]-[X,Y],\ \ \forall\ \ X,Y \in \Gamma
(T^*M).
$$
Then, we have $\frac{\delta}{\delta q^k}t =0,\
\frac{\partial}{\partial p_k}t = g^{0k}$ and $\dot
\nabla_iG_{jk}=0,~\dot \nabla_iH^{jk}= 0$, where
$$
\dot \nabla_iG_{jk}= \frac{\delta}{\delta
q^i}G_{jk}-\Gamma^l_{ij}G_{lk}-\Gamma^l_{ik}G_{lj}
$$
$$
\dot \nabla_iH^{jk}= \frac{\delta}{\delta
q^i}H^{jk}+\Gamma^j_{il}H^{lk}+\Gamma^k_{il}H^{lj}
$$

The above expressions for the components of $N$ can be obtained by
a quite long, straightforward  computation.\\

{\bf Theorem 3.} {\it The almost complex structure $J$ on $T^*M$
is integrable if and only if the base manifold $M$ has constant
sectional curvature $c$ and the function $v$ is given by}
\begin{equation}
v=-\frac{c}{A}.
\end{equation}

{\it Proof.} From the condition $N=0$, one obtains
$$
\{Av(\delta^h_ig_{jk}- \delta^h_jg_{ik})+R^h_{kij}\}p_h=0.
$$

Differentiating with respect to $p_l$, taking $p_h=0~\forall h \in
{\{}1,...,n{\}}$, it follows that the curvature tensor field of
$\dot \nabla$ has the expression

$$
R^l_{kij}=-Av(0)({\delta}^l_ig_{jk}-{\delta}^l_jg_{ik}).
$$
Using by the Schur theorem(in the case when $M$ is connected and
$dim M \geq 3$) it follows that $(M,g)$ has the constant sectional
curvature $c=-Av(0)$ . Then we obtain the expression (11) of $v$.

Conversely, if $(M,g)$ has constant curvature $c$ and $v$ is given
by (11), it
follows in a straightforward way that $N=0$.\\

\bf Remark. \rm  The function $v$ must fulfill the condition
\begin{equation}
A+2tv=\frac{A^2-2ct}{A}>0,~~A>0.
\end{equation}

If $c<0$ then $(T^*M,G,J)$ is a K\"ahler manifold.

If $c>0$ then $(T^*_AM,G,J)$ is a K\"ahler manifold, where
$T^*_AM$ is the tube around zero section in $T^*M$ defined by the
condition $0\leq\|p\|^2<\frac{A^2}{c}.$

 If $c=0$ then $(M,g)$ is a flat
manifold and G becomes a flat metric on $T^*M$.

 The components of the K\"ahler metric $G$ on $T^*M$ are

\begin{equation}
\left\{
\begin{array}{l}
G_{ij}=Ag_{ij}-\frac{c}{A}p_ip_j,
\\ \mbox{ } \\
H^{ij}=\frac{1}{A}g^{ij}+\frac{c}{A(A^2-2ct)}g^{0i}g^{0j}.
\end{array}
\right.
\end{equation}
\newpage
 \vskip5mm {\large \bf 4. A  K\"ahler  Einstein structure on
 $T^*M$}
\vskip5mm
In this section we shall study the property of the
K\"ahler manifold $(T^*M,G,J)$ to be Einstein.

The Levi Civita connection $\nabla$ of the Riemannian manifold
$(T^*M,G)$ is determined by the conditions
$$
\nabla G=0,~~~~~  T =0,
$$
where $T$ is its torsion tensor field. The explicit expression of
this connection is obtained from the formula
$$
2G({\nabla}_XY,Z)=X(G(Y,Z))+Y(G(X,Z))-Z(G(X,Y))+
$$
$$
+G([X,Y],Z)-G([X,Z],Y)-G([Y,Z],X); ~~~~~~ \forall\
X,Y,Z~{\in}~{\Gamma}(T^*M).
$$

The final result can be stated as follows. \\

\bf Theorem 4. {\it The Levi Civita connection ${\nabla}$ of $G$
has the following expression in the local adapted frame
$(\frac{\delta}{\delta q^1},\dots ,\frac{\delta}{\delta
q^n},\frac{\partial}{\partial p_1},\dots ,\frac{\partial}{\partial
p_n}):$
\begin{equation}
\left\{
\begin{array}{l}
 \nabla_\frac{\partial}{\partial
p_i}\frac{\partial}{\partial p_j}
=Q^{ij}_h\frac{\partial}{\partial p_h},\ \ \ \ \ \
\nabla_\frac{\delta}{\delta q^i}\frac{\partial}{\partial
p_j}=-\Gamma^j_{ih}\frac{\partial}{\partial
p_h}+P^{hj}_i\frac{\delta}{\delta q^h},
\\ \mbox{ } \\
\nabla_\frac{\partial}{\partial p_i}\frac{\delta}{\delta
q^j}=P^{hi}_j\frac{\delta}{\delta q^h},\ \ \ \ \ \
\nabla_\frac{\delta}{\delta q^i}\frac{\delta}{\delta
q^j}=\Gamma^h_{ij}\frac{\delta}{\delta
q^h}+S_{hij}\frac{\partial}{\partial p_h},
\end{array}
\right.
\end{equation}
where $Q^{ij}_h, P^{hi}_j, S_{hij}$ are $M$-tensor fields on
$T^*M$, defined by
\begin{equation}
\left\{
\begin{array}{l}
 Q^{ij}_h = \frac{1}{2}G_{hk}(\frac{\partial}{\partial
p_i}H^{jk}+ \frac{\partial}{\partial p_j}H^{ik}
-\frac{\partial}{\partial p_k}H^{ij}),
\\ \mbox{ } \\
 P^{hi}_j=\frac{1}{2}H^{hk}(\frac{\partial}{\partial
p_i}G_{jk}-H^{il}R^0_{ljk}),
\\ \mbox{ } \\
S_{hij}=-\frac{1}{2}G_{hk}\frac{\partial}{\partial
p_k}G_{ij}+\frac{1}{2}R^0_{hij}.
\end{array}
\right.
\end{equation}
\rm

After replacing of the expressions of the involved $M$-tensor
fields , we obtain

\begin{equation}
Q^{ij}_h = \frac{c}{A}H^{ij},\ \ \ \ \ \
P^{hi}_j=-\frac{c}{A}H^{hi},\ \ \ \ \ \
S_{hij}=\frac{c}{A}G_{hj}p_i .
\end{equation}

The curvature tensor field $K$ of the connection $\nabla $ is
obtained from the well known formula
$$
K(X,Y)Z=\nabla_X\nabla_YZ-\nabla_Y\nabla_XZ-\nabla_{[X,Y]}Z,\ \ \
\ \forall\ X,Y,Z\in \Gamma (T^*M).
$$

The components of curvature tensor field $K$ with respect to the
adapted local frame $(\frac{\delta}{\delta q^1},\dots
,\frac{\delta}{\delta q^n},\frac{\partial}{\partial p_1},\dots
,\frac{\partial}{\partial p_n})$ are obtained easily:
\begin{equation}
\left\{
\begin{array}{l}
K(\frac{\delta}{\delta q^i},\frac{\delta}{\delta
q^j})\frac{\delta}{\delta q^k}=\frac{c}{A}(\delta
^h_iG_{jk}-\delta ^h_jG_{ik})\frac{\delta}{\delta q^h},
\\ \mbox{ } \\
K(\frac{\delta}{\delta q^i},\frac{\delta}{\delta
q^j})\frac{\partial}{\partial p_k}=\frac{c}{A}(\delta
^k_jG_{hi}-\delta ^k_iG_{hj})\frac{\partial}{\partial p_h},
\\ \mbox{ } \\
K(\frac{\partial}{\partial p_i},\frac{\partial}{\partial
p_j})\frac{\delta}{\delta q^k}=\frac{c}{A}(\delta
^j_kH^{hi}-\delta ^i_kH^{hj})\frac{\delta}{\delta q^h},
\\ \mbox{ } \\
K(\frac{\partial}{\partial p_i},\frac{\partial}{\partial
p_j})\frac{\partial}{\partial p_k} =\frac{c}{A}(\delta
^i_hH^{jk}-\delta ^j_hH^{ik})\frac{\partial}{\partial p_h},
\\ \mbox{ } \\
K(\frac{\partial}{\partial p_i},\frac{\delta}{\delta
q^j})\frac{\delta}{\delta
q^k}=\frac{c}{A}\delta^i_jG_{hk}\frac{\partial}{\partial p_h},
\\ \mbox{ } \\
K(\frac{\partial}{\partial p_i},\frac{\delta}{\delta
q^j})\frac{\partial}{\partial
p_k}=-\frac{c}{A}\delta^i_jH^{hk}\frac{\delta}{\delta q^h}.
\end{array}
\right.
\end{equation}
{\bf Remark.} From the local coordinates expression of the
curvature tensor field K we obtain that the K\"ahler manifold
$(T^*M,G,J)$ cannot have constant holomorphic sectional curvature.\\

The Ricci tensor field Ric of $\nabla$ is defined by the formula:
$$
Ric(Y,Z)=trace(X\longrightarrow K(X,Y)Z),\ \ \ \forall\  X,Y,Z\in
\Gamma (T^*M).
$$
It follows
$$
\left\{
\begin{array}{l}
Ric(\frac{\delta}{\delta q^i},\frac{\delta}{\delta
q^j})=\frac{cn}{A}G_{ij},
\\ \mbox{ } \\
 Ric(\frac{\partial}{\partial
p_i},\frac{\partial}{\partial p_j})=\frac{cn}{A}H^{ij},
\\ \mbox{ } \\
Ric(\frac{\partial}{\partial p_i},\frac{\delta}{\delta
q^j})=Ric(\frac{\delta}{\delta q^j},\frac{\partial}{\partial
p_i})=0.
\end{array}
\right.
$$
Thus
\begin{equation}
 Ric=\frac{cn}{A}G.
\end{equation}

By using the expressions (17), we have computed the covariant
derivatives of curvature tensor field K in the local adapted frame
$(\frac{\delta}{\delta q^i},\frac{\partial}{\partial p_i})$ with
respect to the connection $\nabla$ and we obtained in all twelve
cases the result is zero.

 Hence we may state our main
result.\\

{\bf Theorem 5.} {\it Assume that the Riemannian manifold $(M,g)$
has constant sectional curvature $c$, the conditions (12) are
fulfilled and the components of the metric $G$ are given by (13).

1. If\ \ $c<0$ then $(T^*M,G,J)$ is a locally symmetric K\"ahler
Einstein manifold.

2. If\ \ $c>0$ then $(T^*_AM,G,J)$ is a locally symmetric K\"ahler
Einstein manifold.}

\vskip 1.5cm

\begin{minipage}{2.5in}
\begin{flushleft}
D.D.Poro\c sniuc\\
Department of Mathematics \\
National College "M. Eminescu" \\
Boto\c sani, Rom\^ania.\\
e-mail: danielporosniuc@lme.ro
\end{flushleft}
\end{minipage}


\begin{thebibliography}{99}

\bibitem{Besse} {\bf Besse, A.,} {\it Einstein manifolds,}~Ergeb.
Math.Grenzgeb.(3) 10, Springer-Verlag, Berlin 1987.
\bibitem{Calabi} {\bf Calabi,~E.,} {\it M\'etriques Kaehl\'eriennes et fibr\'es
holomorphes},~Ann. Scient. Ec. Norm. Sup.,~12 (1979),~269-294.
\bibitem{GheOpr} {\bf Gheorghiev,~Gh.;~Oproiu,~V.,} {\it Variet\u a\c ti
diferen\c tiabile finit \c si infinit dimensionale},~Ed. Academiei
Rom. I (1976),~~II (1979).
\bibitem{KobNom} {\bf Kobayashi, S.; Nomizu, K.,} {\it Foundations
of Differential Geometry I, II,}~Interscience, New York, 1963,
1969.
\bibitem{Lee} {\bf Lee, J.M.,} {\it Ricci. A Mathematica package for doing tensor
calculations in differential geometry. User's Manual, 1992, 2000.}
\bibitem{Mok} {\bf Mok,~K.P.;~Patterson,~E.M.;~Wong,~Y.C.,} {\it Structure
of symmetric tensors of type (0,2) and tensors of type (1,1) on
the tangent bundle},~ Trans. Am. Math. Soc.~234 (1977),~253-278.
\bibitem{Oproiu1} {\bf Oproiu,~V.,} {\it On the differential
geometry of the tangent bundles,} Rev. Roum. Math. Pures Appl. 13
(1968), 847-856.
\bibitem{Oproiu2} {\bf Oproiu,~V.,} {\it A generalization of natural
almost Hermitian structures on the tangent bundles.} Math. J.
Toyama Univ., 22 (1999), 1-14.
\bibitem{Oproiu3} {\bf Oproiu,~V.,} {\it General natural almost Hermitian and
anti-Hermitian structures on the tangent bundles.} Bull. Soc. Sci.
Math. Roum. 43 (91), (2000), 325-340.
\bibitem{Oproiu4} {\bf Oproiu,~V.,} {\it Some new geometric
structures on the tangent bundle}, Public. Math. Debrecen. 55
(1999), 261-281.
\bibitem{Oproiu5} {\bf Oproiu,~V.,} {\it A K\"ahler Einstein
structure on the tangent bundle of a space form}, Int. J. Math.
Math. Sci. 25 (2001), 183-195.
\bibitem{Oproiu6} {\bf Oproiu,~V.,} {\it A locally symmetric
K\"ahler Einstein structure on the tangent bundle of a space
form}, Beitr\"age zur Algebra und Geometrie/ Contributions to
Algebra and Geometry, 40 (1999), 363-372.
\bibitem{OprPap1} {\bf Oproiu,~V.; Papaghiuc,~N.,} {\it A locally
symmetric Kaehler Einstein structure on a tube in the tangent
bundle of a space form}, Revue Roumaine Math. Pures Appl. 45
(2000), 863-871.
\bibitem{OprPap2} {\bf Oproiu,~V.; Papaghiuc,~N.,} {\it A Kaehler structure
on the nonzero tangent bundle of a space form,} Differential Geom.
Appl. 11 (1999), 1-12.
\bibitem{OprPap3} {\bf Oproiu,~V.; Papaghiuc,~N.,} {\it A pseudo-Riemannian
metric on the cotangent bundle,}~An. \c Stiin\c t. Univ. Al. I.
Cuza, Ia\c si 36 (1990), 265-276.
\bibitem{OprPap4} {\bf Oproiu,~V.; Papaghiuc,~N.,} {\it Another pseudo-Riemannian
metric on the cotangent bundle,}~Bull. Inst. Polit., Iasi, 27
(1991), fasc.1-4, sect.1.
\bibitem{OprPor} {\bf Oproiu,~V.;~Poro\c sniuc,~D.D.,} {\it A
K\"ahler Einstein structure on the cotangent bundle of a
Riemannian manifold},~to appear in An. \c Stiin\c t. Univ. Al. I.
Cuza, Ia\c si.
\bibitem{Por1} {\bf Poro\c sniuc,~D.D.,} {\it A locally symmetric
K\"ahler Einstein structure on a tube in the nonzero cotangent
bundle of a space form}, preprint.
\bibitem{Por2} {\bf Poro\c sniuc,~D.D.,} {\it A K\"ahler Einstein structure
 on the nonzero cotangent bundle of a space
form}, preprint.
\bibitem{Sasaki} {\bf Sasaki,~S.,} {\it On the Differential Geometry of
the Tangent Bundle of Riemannian Manifolds},~ Tohoku Math. J.,~10
(1958),~ 238-354.
\bibitem{Tahara2} {\bf Tahara,~M.; Vanhecke,~L.; Watanabe,~Y.,} {\it New
structures on tangent bundles}, Note di Matematica (Lecce), 18
(1998), 131-141.
\bibitem{Tahara} {\bf Tahara,~M.; Watanabe,~Y.,} {\it Natural almost
Hermitian, Hermitian and K\"ahler metrics on the tangent bundles},
Math. J. Toyama Univ., 20 (1997), 149-160.
\bibitem{YanoKob} {\bf Yano,~K.;~Kobayashi,~S.,} {\it Prolongations
of tensor fields and connections to tangent bundles,~ I.~General
Theory}, ~Jour. Math. Soc. Japan,~18 (1966),~194-210.
\bibitem{YanoIsh} {\bf Yano,~K.;~Ishihara,~S.,} {\it Tangent and Cotangent
Bundles},~M. Dekker Inc., New York,~1973.


\end{thebibliography}
\end{document}